\documentclass[preprint]{elsarticle}

\makeatletter
\def\ps@pprintTitle{%
 \let\@oddhead\@empty
 \let\@evenhead\@empty
 \def\@oddfoot{}%
 \let\@evenfoot\@oddfoot}
\makeatother

\usepackage{amsfonts}
\usepackage{amsmath}
\usepackage{amsthm}
\usepackage{amssymb}
\usepackage{lineno,hyperref}
\modulolinenumbers[5]

\usepackage{tikz}
\usetikzlibrary{arrows}
\usetikzlibrary{matrix}

\usepackage[boxed, ruled, vlined]{algorithm2e}

\newcommand{\bo}[1]{\boldsymbol{#1}}

\newcommand{\cov}{\mbox{Cov}}
\newcommand{\E}{\mbox{E}}

\newtheorem{theorem}{Theorem}

\newtheorem{remark}{Remark}

\newtheorem{assumptionvec}{}
\newtheorem{assumptionten}{}

\newcommand\numberthis{\addtocounter{equation}{1}\tag{\theequation}}

\journal{Signal Processing}

\biboptions{authoryear}

\begin{document}

\begin{frontmatter}

\title{Blind source separation of tensor-valued time series\tnoteref{mytitlenote}}
\tnotetext[mytitlenote]{Funding: This work was supported by the Academy of Finland [grant 268703].}

\author{Joni Virta\corref{cor1}}
\ead{joni.virta@utu.fi}
\cortext[cor1]{Corresponding author}

\author{Klaus Nordhausen\corref{cor2}}
\ead{klaus.nordhausen@utu.fi}
\address{Department of Mathematics and Statistics, University of Turku, FI-20014, Finland}

\begin{abstract}
The blind source separation model for multivariate time series generally assumes that the observed series is a linear transformation of an unobserved series with temporally uncorrelated or independent components. Given the observations, the objective is to find a linear transformation that recovers the latent series. Several methods for accomplishing this exist and three particular ones are the classic SOBI and the recently proposed generalized FOBI (gFOBI) and generalized JADE (gJADE), each based on the use of joint lagged moments. In this paper we generalize the methodologies behind these algorithms for tensor-valued time series. We assume that our data consists of a tensor observed at each time point and that the observations are linear transformations of latent tensors we wish to estimate. The tensorial generalizations are shown to have particularly elegant forms and we show that each of them is Fisher consistent and orthogonal equivariant. Comparing the new methods with the original ones in various settings shows that the tensorial extensions are superior to both their vector-valued counterparts and to two existing tensorial dimension reduction methods for i.i.d. data. Finally, applications to fMRI-data and video processing show that the methods are capable of extracting relevant information from noisy high-dimensional data.
\end{abstract}


\end{frontmatter}


\section{Introduction}\label{sec:intro}

\subsection{Blind source separation and time series}

In the classical \textit{blind source separation (BSS)} model one assumes that the observed random vectors $\bo{x}_i$, $i=1,\ldots,n$, are linear transformations of some latent vectors of interest, $\bo{x}_i = \bo{\Omega} \bo{z}_i$, where $\bo{\Omega} \in \mathbb{R}^{p \times p}$ is a full rank \textit{mixing matrix}. Coupling the model with different sets of assumptions on $\bo{z}_i$ gives various well-known models: \textit{(i)} assuming that $\bo{z}_i$ are i.i.d. and have mutually independent components yields the \textit{independent component (IC) model} \citep[see e.g.][]{comon2010handbook}; \textit{(ii)} assuming $\bo{z}_i$ are i.i.d. and have spherical distribution yields an elliptical model for $\bo{x}_i$ \citep{Oja2010} and \textit{(iii)} as a special case of both previous, assuming that $\bo{z}_i$ are i.i.d. and have standard Gaussian distribution yields the general multivariate Gaussian distribution for $\bo{x}_i$.

In the context of time series it is natural to incorporate the time dependency of the data into the model structure and a commonly used BSS model \citep{belouchrani1997blind} assumes that the observed series $\bo{x}_t$ is generated as

\begin{align} \label{eq:v_bssmodel}
\bo{x}_t = \bo{\Omega} \bo{z}_t, \quad t=0, \pm 1, \pm 2,\ldots,
\end{align}
where the latent time series $\bo{z}_t$ satisfies the following three assumptions.

\begin{assumptionvec}\label{assu:v_mean}
$\E[\bo{z}_t] = \bo{0}$
\end{assumptionvec}

\begin{assumptionvec}\label{assu:v_cov}
$\cov[\bo{z}_t] = \bo{I}$
\end{assumptionvec}

\begin{assumptionvec}\label{assu:v_autocov}
$\E[\bo{z}_t \bo{z}_{t+\tau}^T] = \E[\bo{z}_{t+\tau} \bo{z}_{t}^T] = \bo{D}_{\tau}$ is diagonal for all $\tau = 1, 2, \ldots$.
\end{assumptionvec}

Without loss of generality, \ref{assu:v_mean} implies that the observed series is centered and \ref{assu:v_cov} fixes the scales of the columns of $\bo{\Omega}$. After these two assumptions we can still freely change the signs and order of the elements of $\bo{z}_t$ and the corresponding columns of $\bo{\Omega}$ without altering the overall model. Thus the order and the signs of the latent series are unidentifiable which, however, is usually not a problem in practice. \ref{assu:v_mean}--\ref{assu:v_autocov} together also imply that the time series $\bo{z}_t$ and $\bo{x}_t$ are weak second-order stationary and $\bo{x}_t$ satisfies $\E[\bo{x}_t] = \bo{0}$ and $\cov[\bo{x}_t] = \bo{\Omega} \bo{\Omega}^T$.

Time series BSS models such as this have gained in popularity in recent years as general multivariate time series models are demanding in both theory and computation. On the other hand, assuming a BSS model allows the use of well-established univariate time series methods for each of the estimated latent components separately. For some recent contributions see, for example, \citet{OjaKiviluotoMalaroiu2000,WuYu2005,ChenHardleSpokoiny2007,BrodaPaolella2009,LuWuLee2009,ChenHardleSpokoiny2010,GarciaFerrerConzalezPrietoPena2012}.

\subsection{Tensor-valued methods for time series}

In the models discussed above at each time point a $p$-variate vector is observed. Modern data structures are however often more complex. For example, in many applications at each time point data might be observed which is better represented by a tensor. Such applications are for instance spatio-temporal data where at each time point usually a matrix is obtained or fMRI (functional Magnetic Resonance Imaging) data where for each time point a 3-dimensional tensor is recorded. But also video clip data can be seen as a time series where each frame is a matrix- or tensor-valued observation, depending on the number of colors used.

The most common approach to analyzing such data is to convert, following some convention, the tensor into a large vector and then apply standard multivariate methods for vector-valued data. Besides the often practical problem that the vectorized data might be of quite high dimension also information gets lost in this process. As for example \citet{werner2008estimation} point out, after vectorizing the tensors the resulting vectors have a Kronecker structure. Ignoring this structure then means that a much larger number of parameters needs to be estimated.

For i.i.d. data this has recently led to extensive research where methods either model the Kronecker structure or work directly with the tensors. For some recent work on structured multivariate estimation see for example
\citet{werner2008estimation,srivastava2008models,wiesel2012geodesic,greenewald2014robust} and references therein. For some recent contributions for i.i.d. tensor methods see for example
\citet{li2010dimension,zeng-2013,ding2015tensor,zhong-2015} and references therein.

Also independent component analysis (ICA) has already been considered in the context of tensors, some early works being  \citet{beckmann2005tensorial,zhang2008directional,vasilescu2005multilinear}. A fully tensorial model-based approach was however only recently developed in \citet{VirtaLiNordhausenOja2016,VirtaLiNordhausenOja2016b} where tensorial versions of the well-known ICA methods Fourth order blind identification (FOBI) \citep{cardoso1989source} and Joint approximate diagonalization of eigen-matrices (JADE) \citep{cardoso1993blind} were introduced.

Methods for tensor-valued time series seem however not to have had much attention yet although the first steps are for example \citet{WaldenSerroukh2002,RogersLiRussell2013,BahadoriYuLiu2014}. But to the best of our knowledge no tensorial BSS-methods for dependent data have been considered so far. To fill this gap, in this paper we propose tensor extensions for three BSS methods meant for multivariate time series. The first method is called the \textit{Second order blind identification} (SOBI) \citep{belouchrani1997blind} and is based on using second-order information in the form of autocovariance matrices to separate the hidden source series. As such it is best suited for multivariate linear processes and may not work with models having trivial autocovariances, for example with stochastic volatility models such as GARCH. The recently proposed second and third methods, \textit{generalized FOBI} (gFOBI) and \textit{generalized JADE} (gJADE), are the exact opposite and operate on the joint fourth-order moments of the component series, see \citet{MatilainenNordhausenOja2015}. Thus the successful use of either requires non-trivial higher moments, ruling out for example the standard ARMA models. The tensorial extensions of these three methods are respectively called TSOBI, TgFOBI and TgJADE and are discussed in Section \ref{sec:bssten}.

\subsection{Structure of the paper}

In Section \ref{sec:nota} we introduce the used notation and define various concepts of multilinear algebra we need to operate the tensor observations. Although fairly easy to grasp and use, after we define the $m$-flattening of a tensor most tensor operations can be carried out conveniently in a matrix form. Section \ref{sec:bssvec} reviews the theory of SOBI, gFOBI and gJADE and prepares the ground for their tensor versions in Section \ref{sec:bssten} where the corresponding theory and algorithms are discussed. In Section \ref{sec:simu} we first use simulations to compare the presented methods with their vector-valued counterparts for vectorized data and then  use the methods to process simulated fMRI-data and a color video. In both applications the proposed methodology is shown to extract the key elements of the signals in compressed form. In Section \ref{sec:disc} we finally conclude with some discussion and the proofs are gathered in \ref{sec:appe}.

\section{Notation and tensor algebra} \label{sec:nota}

\subsection{Notation in general}

Throughout the paper scalars are denoted by lower-case letters, $a, b, c, \ldots$, vectors by lower-case boldface letters, $\bo{a}, \bo{b}, \bo{c}, \ldots$, matrices by capital boldface letters, $\bo{A}, \bo{B}, \bo{C}, \ldots$, and tensors of general order by capital blackboard letters, $\mathbb{A}, \mathbb{B}, \mathbb{C}, \ldots$ (note that $\mathbb{R}$ still means the real line). The same convention on fonts is followed with random elements, but instead using the letters from the end of the alphabet, $x, y, z, \bo{x}, \bo{y}, \bo{z}$, etc.

\subsection{Regarding matrices}

The standard basis vectors of $\mathbb{R}^p$ are denoted by $\bo{e}_i$, $i=1,\ldots,p$, and using them we can construct the matrices $\bo{E}^{ij} := \bo{e}_i \bo{e}_j^T$, the only non-zero element of $\bo{E}^{ij}$ being the single one as the element $(i, j)$. We further make use of the following sets of $p \times p$ matrices: $\mathcal{P}$, the set of all matrices with a single one in each row and column and rest of the entries zero; $\mathcal{J}$, the set of all diagonal matrices with the diagonal entries equal to $\pm 1$; $\mathcal{D}$, the set of all diagonal matrices with positive diagonal elements and $\mathcal{C}$, the set of all matrices $\bo{C} = \bo{PJD}$ where $\bo{P} \in \mathcal{P}$, $\bo{J} \in \mathcal{J}$ and $\bo{D} \in \mathcal{D}$. The sets $\mathcal{P}, \mathcal{J}$ and $\mathcal{D}$ then respectively correspond to the sets of permutation matrices, heterogeneous sign-change matrices and heterogeneous scaling matrices. Finally, $\| \cdot \|$ is the Frobenius norm and by the equivalence $\bo{A} \equiv \bo{B}$ we mean that $\bo{A} = \bo{PJB}$ for some $\bo{P} \in \mathcal{P}$ and $\bo{J} \in \mathcal{J}$.

\subsection{Regarding tensors}

\begin{figure}
\centering
\includegraphics[width=0.70\textwidth]{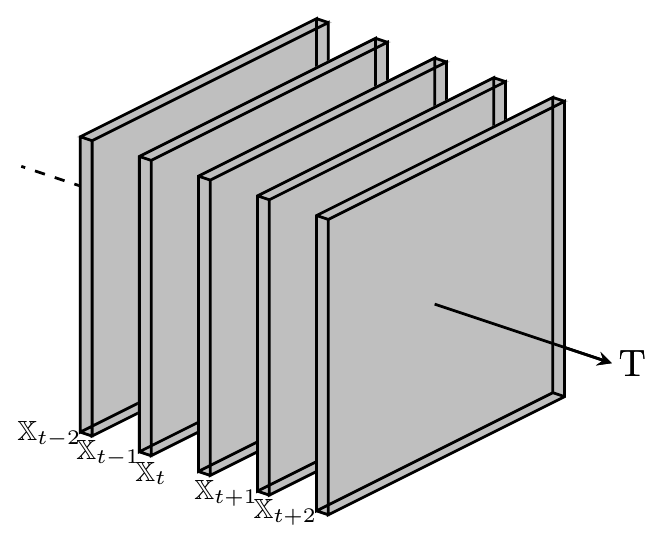}
\caption{Visualization of a tensor-valued time series. In the above scheme a matrix of the same size is observed at each time point and the resulting tensor-valued time series can be thought of as a video with frames corresponding to the individual observed matrices.}
\label{fig:series}
\end{figure}

To manipulate tensors we next provide some basic tools of multilinear algebra, see also \citet{de2000multilinear}. By a tensor-valued time series $\mathbb{X}_t$ we mean the set $\{ \mathbb{X}_t \}_{t=-T}^T$ of realisations of a tensor-valued stochastic process $\mathbb{X}_t \in \mathbb{R}^{p_1 \times \cdots \times p_r}$ on some fixed set of time points $t=-T,\ldots,T$. That is, for each time point we observe a tensor of the same size, something akin to the frames of a video, see Figure \ref{fig:series} for a visual representation.

A useful way of understanding a tensor as a collection of blocks of lower order is provided by the $m$-mode vectors (or fibers), which extend the notion of the columns and rows of a matrix into higher dimensions. The $m$-mode vectors of a tensor $\mathbb{X} \in \mathbb{R}^{p_1 \times \cdots \times p_r}$ are obtained by fixing the values of all other indices but the $m$th and letting the $m$th index vary over its range. The columns and rows of a matrix are thus also its $1$-mode and $2$-mode vectors, respectively. Each $m$-mode vector is of length $p_m$ and the total number of $m$-mode vectors in $\mathbb{X}$ is $\rho_m := \prod_{i \neq m}^r p_i$.  Stacking all of them horizontally in a pre-specified order into a matrix yields what we call the $m$-flattening of a tensor, denoted by $\bo{X}^{(m)} \in \mathbb{R}^{p_m \times \rho_m}$. As is later seen the ordering is not important as long as it is consistent and one option is e.g. the cyclical ordering discussed in \citet{de2000multilinear}.

We next define two multiplication operations, the first one being the $m$th mode linear transformation $(\mathbb{X} \odot_{m} \bo{A}) \in \mathbb{R}^{p_1 \times \cdots \times p_r}$ defined element-wise as
\[\left(\mathbb{X} \odot_{m} \bo{A}\right)_{i_1 \cdots i_r} = \sum_{j_m} x_{i_1 \cdots i_{m-1} j_m i_{m+1} \cdots i_r} a_{i_m j_m},\]
where $\mathbb{X} \in \mathbb{R}^{p_1 \times \cdots \times p_r}$ and $\bo{A} \in \mathbb{R}^{p_m \times p_m}$. The second operation, $(\mathbb{X} \odot_{-m} \mathbb{Y}) \in \mathbb{R}^{p_m \times p_m}$, takes two tensors of the same size, $\mathbb{X}, \mathbb{Y} \in \mathbb{R}^{p_1 \times \cdots \times p_r}$ and acts element-wise as
\[\left(\mathbb{X} \odot_{-m} \mathbb{Y}\right)_{kl} = \sum_{i_1, \ldots ,i_{m-1},i_{m+1}, \ldots ,i_r} x_{i_1 \cdots i_{m-1} k i_{m+1} \cdots i_r} y_{i_1 \cdots i_{m-1} l i_{m+1} \cdots i_r}.\]
Nice interpretations for the previous operations are given by noticing that $(\mathbb{X} \odot_{m} \bo{A})^{(m)} = \bo{A} \bo{X}^{(m)}$ and that $\mathbb{X}\odot_{-m}\mathbb{X} = \bo{X}^{(m)} (\bo{X}^{(m)})^T$. That is, $\mathbb{X} \odot_{m} \bo{A}$ applies the linear transformation given by $\bo{A}$ individually to all $m$-mode vectors of $\mathbb{X}$ and $\mathbb{X}\odot_{-m}\mathbb{X}$ is simply the sum of the outer products of all $m$-mode vectors of $\mathbb{X}$ with themselves, the order we stacked them into $\bo{X}^{(m)}$ clearly not mattering.

Although the above tensor notation is quite straightforward, the familiarity of matrix notation still outweighs its usefulness and we thus choose to work with the $m$-flattened tensors whenever possible. The ability to express everything with matrices also makes the implementation of the algorithms easier.

\section{Blind source separation of multivariate time series} \label{sec:bssvec}

\subsection{The model and the standardized time series}

In this section we review the theory behind three moment-based blind source separation methods for time series, SOBI, gFOBI and gJADE, which are then in Section \ref{sec:bssten} extended for tensor-valued time series.

All three methods assume the model \eqref{eq:v_bssmodel} along with \ref{assu:v_mean} and \ref{assu:v_cov} but each of them has additional assumptions that guarantee the successful identification of the latent time series with the particular method. These extra assumptions will be elaborated more on the next section discussing the estimation of the rotation.

Interestingly, the model and \ref{assu:v_mean} and \ref{assu:v_cov} alone already yield the first step in the estimation of $\bo{z}_t$. This step is given by the whitening of the data, $\bo{x}_t \mapsto \bo{x}_t^{st} := (\bo{\Sigma}_0(\bo{x}_t))^{-1/2} \bo{x}_t$, where $\bo{\Sigma}_0(\bo{x}_t)$ is the covariance matrix of the stationary series $\bo{x}_t$. It is easy to see (e.g. by considering the singular value decomposition of $\bo{\Omega}$) that the whitened series satisfies
\begin{align} \label{eq:v_rotation}
\bo{x}_t^{st} = \bo{U} \bo{z}_t,
\end{align}
for some orthogonal matrix $\bo{U} \in \mathbb{R}^{p \times p}$. The problem is then reduced from that of estimating the $p^2$ parameters of $\bo{\Omega}^{-1}$ into that of estimating the $(1/2)p(p-1)$ parameters of $\bo{U}$, the standardization thus solving half of the problem in terms of number of parameters to estimate.

\subsection{Estimating the rotation}

The unknown rotation $\bo{U}$ is estimated using one of three families of matrix-valued functionals. But prior to defining them we first define another set of functionals that is used in constructing one of the families. $\mathcal{T}$ being a pre-specified set of lags, we define
\begin{align}
\bo{B}_{\tau i j}(\bo{x}_t) &:= \E\left[\bo{e}_i^T \bo{x}_{t+\tau} \bo{x}_{t+\tau}^T \bo{e}_j \cdot \bo{x}_{t}\bo{x}_{t}^T\right], \label{eq:mat_btauij}
\end{align}
where $\tau \in \mathcal{T}$ and $i,j = 1,\ldots,p$. The matrices $\bo{B}_{\tau i j}(\bo{x}^{st}_t)$ then contain specific sums of joint lagged fourth-order moments of the series $\bo{x}^{st}_t$ and, concerning the current estimation problem, unless the underlying series has trivial fourth-order moments they contain information on the unknown rotation.

The three families of interest are then defined as follows.
\begin{align}
\bo{\Sigma}_\tau(\bo{x}_t) &:= \E\left[\bo{x}_t \bo{x}_{t+\tau}^T\right], \label{eq:mat_sigmatau} \\
\bo{B}_\tau(\bo{x}_t) &:= \E\left[\bo{x}_t \bo{x}_{t+\tau}^T \bo{x}_{t+\tau} \bo{x}_{t}^T\right], \label{eq:mat_btau} \\
\bo{C}_{\tau i j}(\bo{x}_t) &:= \bo{B}_{\tau i j}(\bo{x}_t) - \bo{\Sigma}_\tau(\bo{x}_t) (\bo{E}^{ij} + \bo{E}^{ji}) \bo{\Sigma}_\tau(\bo{x}_t)^T - \delta_{ij} \textbf{I}, \label{eq:mat_ctauij}
\end{align}
where $\tau \in \mathcal{T}$, $i,j = 1,\ldots,p$ and $\delta_{ij}$ is the Kronecker delta. Consider then the following three additional assumptions.

\begin{assumptionvec}\label{assu:v_sobiworks}
The matrix $\bo{\Sigma}_\tau(\bo{z}_t)$ is diagonal for all $\tau \in \mathcal{T}$.
\end{assumptionvec}

\begin{assumptionvec}\label{assu:v_gfobiworks}
The matrix $\bo{B}_\tau(\bo{z}_t)$ is diagonal for all $\tau \in \mathcal{T}$.
\end{assumptionvec}

\begin{assumptionvec}\label{assu:v_gjadeworks}
The matrix $\bo{C}_{\tau i j}(\bo{z}_t)$ is diagonal for all $\tau \in \mathcal{T}$, $i,j =1,\ldots,p$.
\end{assumptionvec}

If we then assume the model \eqref{eq:v_bssmodel}, \ref{assu:v_mean} and \ref{assu:v_cov} and additionally either \ref{assu:v_sobiworks}, \ref{assu:v_gfobiworks} or \ref{assu:v_gjadeworks} it can be shown that the matrix $\bo{U}^T$ diagonalizes respectively the matrices $\bo{\Sigma}_\tau(\bo{x}^{st}_t)$, $\bo{B}_\tau(\bo{x}^{st}_t)$ or $\bo{C}_{\tau i j}(\bo{x}^{st}_t)$ for all $\tau \in \mathcal{T}$, $i,j=1,\ldots,p$. Technically the eigendecomposition of any one of them could then under the suitable assumptions be used to estimate the needed rotation. However, for the rotation to be unique the eigenvalues of the considered matrix need to be distinct and this is not necessarily the case. One solution is, instead of diagonalizing just one of the matrices, to approximately diagonalize several of them using the so-called joint approximate diagonalization. A often used technique for performing the joint diagonalization of a chosen set of matrices $\mathcal{M}$ is given by the optimization problem
\begin{align} \label{eq:simul_diag}
\underset{\bo{U}\bo{U}^T = \bo{I}}{\mbox{max}}\sum_{\bo{M} \in \mathcal{M}} \| \mbox{diag}(\bo{U}^T \bo{M} \bo{U}) \|^2,
\end{align}
that is, we find the orthogonal matrix that makes all matrices in the set as diagonal as possible, as measured by the sum of squares of the resulting diagonal values. It can be shown that this is equivalent to minimizing the sum of squared off-diagonal values. Choosing $\mathcal{M} = \{ \bo{\Sigma}_\tau (\bo{x}^{st}_t) \}_{\tau \in \mathcal{T}}$, $\mathcal{M} = \{ \bo{B}_\tau (\bo{x}^{st}_t) \}_{\tau \in \mathcal{T}}$ or $\mathcal{M} = \{ \bo{C}_{\tau ij} (\bo{x}^{st}_t) \}_{\tau \in \mathcal{T}, i,j=1,\ldots,p}$ gives then SOBI, gFOBI and gJADE, respectively. For further details on the joint diagonalization, see \citet{belouchrani1997blind,IllnerMiettinenFuchsTaskienNordhausenOjaTheis2015}.

Intuitively, diagonalizing a set of matrices means that no single matrix needs to carry the information to separate all $p$ components as long as the information exists in at least one of the matrices in the set. This reasoning is encoded in the following three assumptions that guarantee that the correct solution is actually found in SOBI, gFOBI and gJADE, respectively.

\begin{assumptionvec}\label{assu:v_sobiunique}
For all pairs $i \neq j$ there exists $\tau \in \mathcal{T}$ such that the $i$th and $j$th diagonal elements of $\bo{\Sigma}_\tau (\bo{z}_t)$ are distinct.
\end{assumptionvec}

\begin{assumptionvec}\label{assu:v_gfobiunique}
For all pairs $i \neq j$ there exists $\tau \in \mathcal{T}$ such that the $i$th and $j$th diagonal elements of $\bo{B}_\tau (\bo{z}_t)$ are distinct.
\end{assumptionvec}

\begin{assumptionvec}\label{assu:v_gjadeunique}
For $p-1$ components $i$ there exists $\tau \in \mathcal{T}$ such that $\bo{C}_{\tau ii} (\bo{z}_t) \neq \bo{0}$.
\end{assumptionvec}

The nature of the above assumptions shows that SOBI is targeted to a different set of distributions for $\bo{z}$ than gFOBI and gJADE; SOBI uses second order information in the form of autocovariances to separate the latent time series while gFOBI and gJADE implicitly assume no second order information exists and move directly to the use fourth-order moments.

To put everything together, all algorithms require the model \eqref{eq:v_bssmodel} along with \ref{assu:v_mean} and \ref{assu:v_cov} and additionally SOBI assumes \ref{assu:v_sobiworks} and \ref{assu:v_sobiunique}; gFOBI assumes \ref{assu:v_gfobiworks} and \ref{assu:v_gfobiunique} and gJADE assumes \ref{assu:v_gjadeworks} and \ref{assu:v_gjadeunique}. Note that the standard assumption \ref{assu:v_autocov} implies the weaker one \ref{assu:v_sobiworks} and can replace it if needed. Also, since \ref{assu:v_gfobiworks} and \ref{assu:v_gjadeworks} have little meaning in applications sometimes it might be more practical to assume the following stronger condition that implies all of \ref{assu:v_autocov}, \ref{assu:v_sobiworks}, \ref{assu:v_gfobiworks} and \ref{assu:v_gjadeworks}.

\begin{assumptionvec}\label{assu:v_indep}
The component series of $\bo{z}_t$ are independent.
\end{assumptionvec}

Having estimated the rotation with one of the methods we define the related \textit{unmixing matrix functional} as
\begin{align}
\bo{\Gamma}(\bo{x}_t) := \bo{U}^T(\bo{x}_t) (\bo{\Sigma}_0(\bo{x}_t))^{-1/2},
\end{align}
where $\bo{U}^T(\bo{x}_t)$ is the diagonalizing rotation estimated from the series $\bo{x}_t$ notated as a functional of the data. The estimated components are then obtained via the transformation $\bo{x}_t \mapsto \bo{\Gamma}(\bo{x}_t) \bo{x}_t$. One important property of SOBI, gFOBI and gJADE is that they are all affine equivariant, that is, the related unmixing matrix functionals satisfy $\bo{\Gamma}(\bo{A} \bo{x}_t) \equiv \bo{\Gamma}(\bo{x}_t) \bo{A}^{-1}$ for all non-singular $\bo{A}$, even for $\bo{x}_t$ outside the model \eqref{eq:v_bssmodel}. Thus the choice of the coordinate system does not affect the resulting estimated component series and the three functionals can be seen as examples of invariant coordinate system (ICS) functionals, see \citet{TylerCritchleyDumbgenOja:2009}.

\section{Tensor blind source separation} \label{sec:bssten}

\subsection{The model and the standardized time series}

Using the linear transformation operator $\odot_m$ we define the  tensor extension of the model \eqref{eq:v_bssmodel} in a manner similar to the ICA model extensions in \citet{VirtaLiNordhausenOja2016,VirtaLiNordhausenOja2016b}. This yields the following \textit{tensor blind source separation model} for the tensor-valued time series $\mathbb{X}_t$.
\begin{align} \label{eq:t_bssmodel}
\mathbb{X}_t = \mathbb{Z}_t \odot_1 \bo{\Omega}_1 \cdots \odot_r \bo{\Omega}_r, \quad t = 0, \pm1, \pm2,\ldots,
\end{align}
where $\bo{\Omega}_m \in \mathbb{R}^{p_m \times p_m}$, $m=1,\ldots,r$, are non-singular $m$-mode mixing matrices and the unobserved tensor-valued time series $\mathbb{Z}_t$ satisfies the following three assumptions for all $t = 0, \pm1, \pm2,\ldots$

\begin{assumptionten}\label{assu:t_mean}
$\E[vec(\mathbb{Z}_t)] = \bo{0}$.
\end{assumptionten}

\begin{assumptionten}\label{assu:t_cov}
$\cov[vec(\mathbb{Z}_t)] = \bo{I}$.
\end{assumptionten}

\begin{assumptionten}\label{assu:t_autocov}
$\E[\mathbb{Z}_t^{(m)} (\mathbb{Z}_{t+\tau}^{(m)})^T] = \E[\mathbb{Z}_{t+\tau}^{(m)} (\mathbb{Z}_{t}^{(m)})^T] = \bo{D}^m_{\tau}$ are diagonal for all $m=1,\ldots,r$ and for all $\tau = 1,2,\ldots$
\end{assumptionten}

Assumption \ref{assu:t_mean} again without loss of generality requires that the observed time series is centered and Assumption \ref{assu:t_cov} fixes the relative scales of the columns of the mixing matrices. Note that one can still multiply any of the mixing matrices by a constant and divide one of them by the same constant without changing the model. Assumption \ref{assu:t_cov} further implies that the covariance matrix of the vectorized series $vec(\mathbb{X}_t)$ has the Kronecker structure, $\cov(vec(\mathbb{X}_t)) = (\bo{\Omega}_m \bo{\Omega}_m^T) \otimes \cdots \otimes (\bo{\Omega}_1 \bo{\Omega}_1^T)$, where $\otimes$ is the Kronecker product. Assumption \ref{assu:t_autocov} is a straight generalization of \ref{assu:v_autocov} and, recalling the definition of $\bo{Z}^{(m)}_{t}$, in principle says that the $m$-mode vectors of $\mathbb{Z}_t$ are ``on average stationary''. Like for the vector methods, also now we need additional assumptions to guarantee that each of the methods finds a solution and that the solution is the correct one and these are discussed in the next section.

Curiously, the BSS model \eqref{eq:t_bssmodel} can be viewed as a collection of low-multilinear rank models for each observed tensor $\mathbb{Z}_t$ with the constraint that the factor matrices are equal for each time point. Still another viewpoint is to treat model \eqref{eq:t_bssmodel} as a low-multilinear rank model for the whole tensor-valued time series viewed as a single tensor of order $r + 1$, with the time dimension left uncompressed. If we further impose the additional constraint that the mixing matrices be orthogonal the model can be seen as a multilinear singular value decomposition in the same two different ways.

Two interesting analogies between the model \eqref{eq:t_bssmodel} and \textit{tensor segmentation} can also be drawn. In \cite{bousse2017tensor} it is assumed that in vector-valued BSS either the individual signals or the columns of the mixing matrix (and later both) admit low-rank tensorizations. Consider first the former case and assume that we have a data matrix $\textbf{X} \in \mathbb{R}^{T \times p}$ with the $p$ columns corresponding to the $p$ individual observed signals. Then our working assumption is that the \textit{rows} of $\textbf{X}$ are actually vectorized tensors while the segmentation conversely assumes that the \textit{columns} of $\textbf{X}$ are actually vectorized tensors. Consider next the case where the columns of the mixing matrix are actually vectorized tensors and therein a subcase where each of the tensors is of rank one. Then the vector BSS-model \eqref{eq:v_bssmodel} can be written as
\begin{align}\label{eq:khatrirao}
\textbf{x}_t = (\boldsymbol{\Omega}^*_r \circledcirc \cdots \circledcirc \boldsymbol{\Omega}^*_1) \textbf{z}_t,
\end{align}
where the vectors and matrices are of suitable size, the $j$th columns of the matrices $\boldsymbol{\Omega}^*_m$, $m = 1, \ldots ,r$, give the factors of the rank-1 tensorization of the $j$th column of the mixing matrix $\boldsymbol{\Omega}$, $j = 1, \ldots, p$ and $\circledcirc$ denotes the Khatri-Rao product (column-wise Kronecker product). Comparing \eqref{eq:khatrirao} with the vectorization of the model \eqref{eq:t_bssmodel},
\[vec(\mathbb{X}_t) = (\boldsymbol{\Omega}_r \otimes \cdots \otimes \boldsymbol{\Omega}_1) vec(\mathbb{Z}_t),\]
shows that the two approaches are linked, both assuming a certain structure on the mixing matrix.

As in the vector case, the first step towards solving \eqref{eq:t_bssmodel} is given by a standardization based on second-order moments. Define first the $m$th mode covariance matrices as
\begin{align}
\bo{\Sigma}_0^m(\mathbb{X}_t) := \frac{1}{\rho_m}\E\left[\bo{X}^{(m)}_{t} (\bo{X}^{(m)}_{t})^T \right],
\end{align}
where $m=1,\ldots,r$. The matrix $\bo{\Sigma}_0^m(\mathbb{X}_t)$ is essentially the mean covariance matrix of all $m$-mode vectors of $\mathbb{X}_t$ and thus measures the average spatial dependency in the $m$th direction. The first step of the unmixing procedure is then given by a simultaneous standardization from all modes of $\mathbb{X}_t$:
\begin{align}
\mathbb{X}^{st}_t := \mathbb{X}_t \odot_1 \left(\bo{\Sigma}_0^1(\mathbb{X}_t)\right)^{-1/2} \cdots \odot_r \left(\bo{\Sigma}_0^r(\mathbb{X}_t)\right)^{-1/2},
\end{align}
where $\left(\bo{\Sigma}_0^m(\mathbb{X}_t)\right)^{-1/2}$ is the unique symmetric inverse square root of $\bo{\Sigma}_0^m(\mathbb{X}_t)$, $m=1, \ldots ,r$. This particular standardization was first considered by \citet{VirtaLiNordhausenOja2016} in the context of independent component analysis and there the authors also gave the following theorem (in the guise of the ICA model).
\begin{theorem} \label{theo:t_rotation}
Assume that the tensor-valued time series $\mathbb{X}_t$ is generated by the model \eqref{eq:t_bssmodel} and satisfies Assumptions \ref{assu:t_mean} and \ref{assu:t_cov}. Then
\begin{align}
\mathbb{X}^{st}_t \propto \mathbb{Z}_t \odot_1 \bo{U}_1 \cdots \odot_r \bo{U}_r,
\end{align}
for some orthogonal matrices $\bo{U}_1 \in \mathbb{R}^{p_1 \times p_1}$, \ldots, $\bo{U}_r \in \mathbb{R}^{p_r \times p_r}$. The constant of proportionality is $(\prod_{m=1}^r p_m^{1/2})^{r-1} \| \bo{\Omega}_r \otimes \cdots \otimes \bo{\Omega}_1 \|^{1-r}$.
\end{theorem}
The proof of the above theorem does not depend in any way on the time series nature of $\mathbb{X}_t$ and is essentially similar to that of Theorem 5.3.1 in \citet{VirtaLiNordhausenOja2016} and thus omitted here. Comparing the result of Theorem \ref{theo:t_rotation} to \eqref{eq:v_rotation} shows that the former serves as a direct tensor analogy to the latter and again reduces the problem of inverting all matrices $\bo{\Omega}_m$ into the problem of inverting the unknown rotations $\bo{U}_m$. 

\subsection{Estimating the rotations}

For estimating the unknown rotations we use tensor counterparts for the families of functionals \eqref{eq:mat_sigmatau}, \eqref{eq:mat_btau} and \eqref{eq:mat_ctauij}, but first we have to define a counterpart for \eqref{eq:mat_btauij}. Assuming a fixed mode $m$ and fixed lags $\tau_1$, $\tau_2$, $\tau_3$ and $\tau_4$ we define
\begin{align}
\bo{B}_{\tau_1 \tau_2 \tau_3 \tau_4 i j}^m (\mathbb{X}_t) &:= \frac{1}{\rho_m}\E\left[\left(\bo{e}_i^T \bo{X}^{(m)}_{t+\tau_1}\left( \bo{X}^{(m)}_{t+\tau_2} \right)^T \bo{e}_j \right) \cdot \bo{X}^{(m)}_{t+\tau_3}\left( \bo{X}^{(m)}_{t+\tau_4} \right)^T\right],
\end{align}
where $i,j = 1,\ldots,p_m$. Using the above we can then define the functionals that allow the estimation of the $r$ rotations. Each of the following is a $p_m \times p_m$ matrix-valued functional and we again assume that $\mathcal{T}$ is a fixed set of lags.

\begin{align*}
\bo{\Sigma}_\tau^m (\mathbb{X}_t) :=& \frac{1}{\rho_m}\E\left[\bo{X}^{(m)}_t\left( \bo{X}^{(m)}_{t+\tau} \right)^T\right], \numberthis \label{eq:mat_sigmamtau} \\
\bo{B}_\tau^m (\mathbb{X}_t) :=& \frac{1}{\rho_m}\E\left[\bo{X}^{(m)}_t\left( \bo{X}^{(m)}_{t+\tau} \right)^T \bo{X}^{(m)}_{t+\tau}\left( \bo{X}^{(m)}_{t} \right)^T\right], \numberthis \label{eq:mat_bmtau} \\
\bo{C}_{\tau i j}^m (\mathbb{X}_t) :=& \bo{B}_{0 \tau \tau 0 i j}^m (\mathbb{X}_t) + \bo{B}_{0 \tau 0 \tau i j}^m (\mathbb{X}_t) - \bo{B}_{\tau \tau 0 0 i j}^m (\mathbb{X}_t) \numberthis \label{eq:mat_cmtauij} \\
&- \bo{\Sigma}_0^m (\mathbb{X}_t)(\bo{E}^{ij} + \bo{E}^{ji} + \bo{I})\bo{\Sigma}_0^m (\mathbb{X}_t)^T,
\end{align*}
where $\tau \in \mathcal{T}$ and $i,j = 1,\ldots,p_m$. The matrices $\bo{B}_\tau^m (\mathbb{X}_t)$ and $\bo{C}_{\tau i j}^m (\mathbb{X}_t)$ can be seen as lagged versions of the matrices used in TFOBI and TJADE, see \citet{VirtaLiNordhausenOja2016,VirtaLiNordhausenOja2016b}. Note that the basic form of the matrices $\bo{C}_{\tau i j}^m (\mathbb{X}_t)$ in \eqref{eq:mat_cmtauij} differs from the form of their vector counterparts in \eqref{eq:mat_ctauij}. This modification is needed to obtain the results in Theorem \ref{theo:t_diag} below.

We further introduce the following three assumptions, relating to the matrices \eqref{eq:mat_sigmamtau}, \eqref{eq:mat_bmtau} and \eqref{eq:mat_cmtauij}, respectively.

\begin{assumptionten}\label{assu:t_tsobiworks}
For all $m=1, \ldots ,r$ and $\tau \in \mathcal{T}$ the matrix $\bo{\Sigma}^m_\tau(\mathbb{Z}_t)$ is diagonal.
\end{assumptionten}

\begin{assumptionten}\label{assu:t_tgfobiworks}
For all $m=1, \ldots ,r$ and $\tau \in \mathcal{T}$ the matrix $\bo{B}^m_\tau(\mathbb{Z}_t)$ is diagonal.
\end{assumptionten}

\begin{assumptionten}\label{assu:t_tgjadeworks}
For all $m=1, \ldots ,r$, $i,j =1,\ldots,p_m$ and $\tau \in \mathcal{T}$ the matrix $\bo{C}^m_{\tau i j}(\mathbb{Z}_t)$ is diagonal.
\end{assumptionten}

Putting all the previous concepts together, a connection between the functionals \eqref{eq:mat_sigmamtau}, \eqref{eq:mat_bmtau} and \eqref{eq:mat_cmtauij} and the estimation of the rotations is stated in the next theorem.

\begin{theorem} \label{theo:t_diag}
Assume that the tensor-valued time series $\mathbb{X}_t$ is generated by the model \eqref{eq:t_bssmodel} and satisfies Assumptions \ref{assu:t_mean}, \ref{assu:t_cov} and additionally either $i)$ \ref{assu:t_tsobiworks}, $ii)$ \ref{assu:t_tgfobiworks} or $iii)$ \ref{assu:t_tgjadeworks}. Then, for a fixed mode $m=1,\ldots,r$, the matrix $\bo{U}_m^T$, where $\bo{U}_m$ is the $m$th true rotation in Theorem \ref{theo:t_rotation}, diagonalizes respectively $i)$ the matrices $\bo{\Sigma}_\tau^m (\mathbb{X}^{st}_t)$, $\tau \in \mathcal{T}$, $ii)$ the matrices $\bo{B}_\tau^m (\mathbb{X}^{st}_t)$, $\tau \in \mathcal{T}$ or $iii)$ the matrices $\bo{C}_{\tau i j}^m (\mathbb{X}^{st}_t)$, $\tau \in \mathcal{T}$, $i,j=1,\ldots,p_m$.
\end{theorem}

Following the ideas of Section~\ref{sec:bssvec}, we can then estimate the $m$th rotation $\bo{U}_m^T$ by simultaneously diagonalizing, see (\ref{eq:simul_diag}), one of the matrix sets: $\mathcal{M}_m = \{ \bo{\Sigma}_\tau^m (\mathbb{X}^{st}_t) \}_{\tau \in \mathcal{T}}$, $\mathcal{M}_m = \{ \bo{B}_\tau^m (\mathbb{X}^{st}_t) \}_{\tau \in \mathcal{T}}$ or $\mathcal{M}_m = \{ \bo{C}^m_{\tau ij} (\mathbb{X}^{st}_t) \}_{\tau \in \mathcal{T}, i,j=1,\ldots,p}$, yielding TSOBI, TgFOBI and TgJADE, respectively. Note however that Theorem \ref{theo:t_diag} does not guarantee that the estimated diagonalizing rotation is necessarily the correct one (it might be a different basis for the correct eigenspace) and to ensure that we further need the following assumptions which apply to TSOBI, TgFOBI and TgJADE, respectively.

\begin{assumptionten}\label{assu:t_tsobiunique}
For all $m=1,\ldots,r$ and for all pairs $i \neq j$ there exists $\tau \in \mathcal{T}$ such that the $i$th and $j$th diagonal elements of $\bo{\Sigma}^m_\tau(\mathbb{Z}_t)$ are distinct.
\end{assumptionten}

\begin{assumptionten}\label{assu:t_tgfobiunique}
For all $m=1,\ldots,r$ and for all pairs $i \neq j$ there exists $\tau \in \mathcal{T}$ such that the $i$th and $j$th diagonal elements of $\bo{B}^m_\tau(\mathbb{Z}_t)$ are distinct.
\end{assumptionten}

\begin{assumptionten}\label{assu:t_tgjadeunique}
For all $m=1,\ldots,r$, for $p_m-1$ indices $i$ there exists $\tau \in \mathcal{T}$ such that $\bo{C}^m_{\tau ii} (\mathbb{Z}_t) \neq \bo{0}$.
\end{assumptionten}

Putting the assumptions again together all three methods assume the model \eqref{eq:t_bssmodel} with \ref{assu:t_mean} and \ref{assu:t_cov}. Additionally, TSOBI requires \ref{assu:t_tsobiworks} and \ref{assu:t_tsobiunique}; TgFOBI requires \ref{assu:t_tgfobiworks} and \ref{assu:t_tgfobiunique} and TgJADE requires \ref{assu:t_tgjadeworks} and \ref{assu:t_tgjadeunique}. Since the assumptions \ref{assu:t_tsobiworks}, \ref{assu:t_tgfobiworks} and \ref{assu:t_tgjadeworks} regarding the working of the methods are again quite impractical it is worth noting that the general assumption \ref{assu:t_autocov} again implies \ref{assu:t_tsobiworks} and all of \ref{assu:t_autocov}, \ref{assu:t_tsobiworks}, \ref{assu:t_tgfobiworks} and \ref{assu:t_tgjadeworks} are implied by the following stronger, but more intuitive assumption.

\begin{assumptionten}\label{assu:t_indep}
The component series of $\mathbb{Z}_t$ are independent.
\end{assumptionten}

This implication can be proven by inspecting the matrices $\bo{\Sigma}^m_\tau(\mathbb{Z}_t)$, $\bo{B}^m_\tau(\mathbb{Z}_t)$ and $\bo{C}^m_{\tau ij} (\mathbb{Z}_t)$ element-wise under \ref{assu:t_indep}, e.g. in the manner of proof of Theorem 1 in \citet{VirtaLiNordhausenOja2016b}.

An algorithm for applying the three proposed methods is given next. The choices for the set of matrices $\mathcal{M}_m$ leading to the different methods are listed above after Theorem \ref{theo:t_diag}.

\begin{algorithm}[H]
\SetAlgoLined
Center: $\mathbb{X}_t \leftarrow \mathbb{X}_t - \bar{\mathbb{X}}_t$\;
Compute $\bo{\Sigma}^m_0(\mathbb{X}_t)$, $m=1,\ldots,r$\;
Standardize: $\mathbb{X}_t \leftarrow \mathbb{X}_t \odot_1 (\bo{\Sigma}^1_0(\mathbb{X}_t))^{-1/2} \cdots \odot_r (\bo{\Sigma}^r_0(\mathbb{X}_t))^{-1/2}$\;
Choose the sets $\mathcal{M}_m$\;
For all $m=1,\ldots,r$, jointly diagonalize $\mathcal{M}_m$ obtaining $\bo{U}^T_m$\;
Rotate: $\mathbb{X}_t \leftarrow \mathbb{X}_t \odot_1 \bo{U}^T_1 \cdots \odot_r \bo{U}^T_r$\;
\caption{TSOBI, TgFOBI and TgJADE. \label{algo:algo1}}
\end{algorithm}

\begin{remark}
Like in the original SOBI, the matrices $\bo{\Sigma}^m_\tau(\mathbb{X}^{st}_t)$ are in practice rarely symmetric, thus not satisfying \ref{assu:t_tsobiworks}, and in the algorithm one should instead use the symmetrized versions, $(1/2)(\bo{\Sigma}^m_\tau(\mathbb{X}^{st}_t) + \bo{\Sigma}^m_\tau(\mathbb{X}^{st}_t)^T)$. Whichever we use makes theoretically no difference as both are under the model diagonalized by the same orthogonal matrix. This practice is described also for other SOBI variants e.g. in \citet{miettinen2014deflation,ilmonen2015affine}.
\end{remark}

As in the vector case, we can for tensor-valued time series define the $m$-mode unmixing functionals as
\begin{align}
\bo{\Gamma}^m(\mathbb{X}_t) = \bo{U}_m^T(\mathbb{X}_t) (\bo{\Sigma}^m_0(\mathbb{X}_t))^{-1/2},
\end{align}
where $\bo{U}_m^T(\mathbb{X}_t)$ is the diagonalizing $m$th rotation estimated from the series $\mathbb{X}_t$ notated as a functional of the data. Using the unmixing functionals of all $r$ modes the final solution is then obtained via the transformation $\mathbb{X}_t \mapsto \mathbb{X}_t \odot_1 \bo{\Gamma}^1(\mathbb{X}_t) \cdots \odot_r \bo{\Gamma}^r(\mathbb{X}_t)$. Of the tensor extensions, however, none is affine equivariant in the same sense as their vector counterparts if $r > 1$. That is, it is not in general true that $\bo{\Gamma}^m(\mathbb{X}_t \odot_1 \bo{A}_1 \cdots \odot_r \bo{A}_r) \equiv \bo{\Gamma}^m(\mathbb{X}_t) \bo{A}_m^{-1}$ for all non-singular $\bo{A}_1, \ldots, \bo{A}_r$. However, the previous holds when all $\bo{A}_1, \ldots, \bo{A}_r$ are orthogonal and all the tensorial extensions are thus \textit{orthogonally equivariant}. See \citet{VirtaLiNordhausenOja2016} for more discussion and a conjecture that affine equivariance is unobtainable in the general tensor case.

\begin{remark}
Using no temporal information, that is, choosing $\mathcal{T} = \{0\}$ in either TgFOBI or TgJADE yields respectively the tensorial ICA methods TFOBI and TJADE as introduced in \citet{VirtaLiNordhausenOja2016,VirtaLiNordhausenOja2016b}.
\end{remark}

\begin{remark}
In \citet{VirtaLiNordhausenOja2016,VirtaLiNordhausenOja2016b} it is shown that both TFOBI and TJADE have alternative versions for \eqref{eq:mat_bmtau} and \eqref{eq:mat_cmtauij} and such modifications could easily be formulated also for TgFOBI and TgJADE. However, \citet{VirtaLiNordhausenOja2016} also showed using asymptotic properties that for the alternative form of TFOBI to have lower total asymptotic variance than the main form the majority of the components of the latent tensor need to have nearly symmetric distributions. \citet{VirtaLiNordhausenOja2016b} on the other hand showed that both versions of TJADE produce almost identical results, the alternative form being computationally much slower. Based on these considerations we thus choose to work with the primary forms in \eqref{eq:mat_bmtau} and \eqref{eq:mat_cmtauij}.
\end{remark}

\section{Examples} \label{sec:simu}

\subsection{Simulation study}

We next compared the performances of the proposed methods to those of their vector-valued counterparts by applying the former to the actual simulated tensor-valued data and the latter to the same data in vectorized form. The proposed methods take advantage of the tensor structure and are thus expected to perform better. The simulations were done with R \citep{Rcore} using the packages \textit{fGarch} \citep{RfGarch}, \textit{ggplot2} \citep{ggplot2}, \textit{JADE} \citep{MiettinenNordhausenTaskinen2017}, \textit{stochvol} \citep{Rstochvol} and \textit{tsBSS} \citep{RtsBSS}. Additionally, our package \textit{tensorBSS} \citep{RtensorBSS} provides implementations of all the discussed tensor methods.

We considered two settings both consisting of simulated tensor-valued time series $\mathbb{Z}_t \in \mathbb{R}^{3 \times 2 \times 2}$ with the lengths of either $T = 1000$, $2000$, $4000$, $8000$, $16000$ or $32000$. In the first setting all component series of $\mathbb{Z}_t$ were generated from different ARMA models and in the second setting from various models exhibiting stochastic volatility (SV), see \ref{sec:appe} for detailed information. The tensors $\mathbb{Z}_t$ were then mixed as $\mathbb{X}_t := \mathbb{Z}_t \odot_1 \bo{A}_1 \odot_2 \bo{A}_2 \odot_3 \bo{A}_3$ where the square matrices $\bo{A}_1$, $\bo{A}_2$ and $\bo{A}_3$ all either had elements sampled independently from the standard normal distribution or were random orthogonal matrices uniform with respect to the Haar measure.

In addition to the proposed methods also the non-lagged versions of TgFOBI and TgJADE, namely TFOBI and TJADE, along with their vector counterparts FOBI and JADE were included in the simulation to see whether the lag information contributes anything to the estimation. In total we thus had 10 methods, FOBI, JADE, SOBI, gFOBI and gJADE and the tensorial version of each. Furthermore, for the SOBI-based methods we used the lag set $\mathcal{T} = \{1,2, \ldots ,12\}$ and for methods based on gFOBI and gJADE the lag set $\mathcal{T} = \{0,1, \ldots ,12\}$, see \citet{MatilainenNordhausenOja2015} for previous use of these particular lag sets. For the reasons discussed in the introductory section we expect TSOBI to outperform the other methods in the ARMA case and TgJADE to do the same in the stochastic volatility case.

To make the results of the vector-valued and tensor-valued methods comparable Kronecker products, $\hat{\bo{\Gamma}}{}^3 \otimes \hat{\bo{\Gamma}}{}^2 \otimes \hat{\bo{\Gamma}}{}^1$, of the 3 unmixing matrix functionals produced by each of the tensor methods were taken. Thus both the Kronecker products and the matrices estimated by the vector methods estimate the inverse of the same $12 \times 12$ matrix $\bo{A}_3 \otimes \bo{A}_2 \otimes \bo{A}_1$. We then measured how close each estimate is to the true inverse via the minimum distance index \citep{ilmonen2010new}:
\begin{align*}
MDI(\hat{\bo{\Gamma}}) = \frac{1}{\sqrt{p - 1}} \underset{\bo{C} \in \mathcal{C}}{\text{inf}}\| \bo{C} \hat{\bo{\Gamma}} (\bo{A}_3 \otimes \bo{A}_2 \otimes \bo{A}_1) - \bo{I} \|,
\end{align*}
where $p$ is now 12 and $\hat{\bo{\Gamma}}$ is the estimated Kronecker product or unmixing matrix. The value of the index can be shown to vary from 0 to 1, the former indicating perfect separation.

\begin{figure}
\centering
\includegraphics[width=1.00\textwidth]{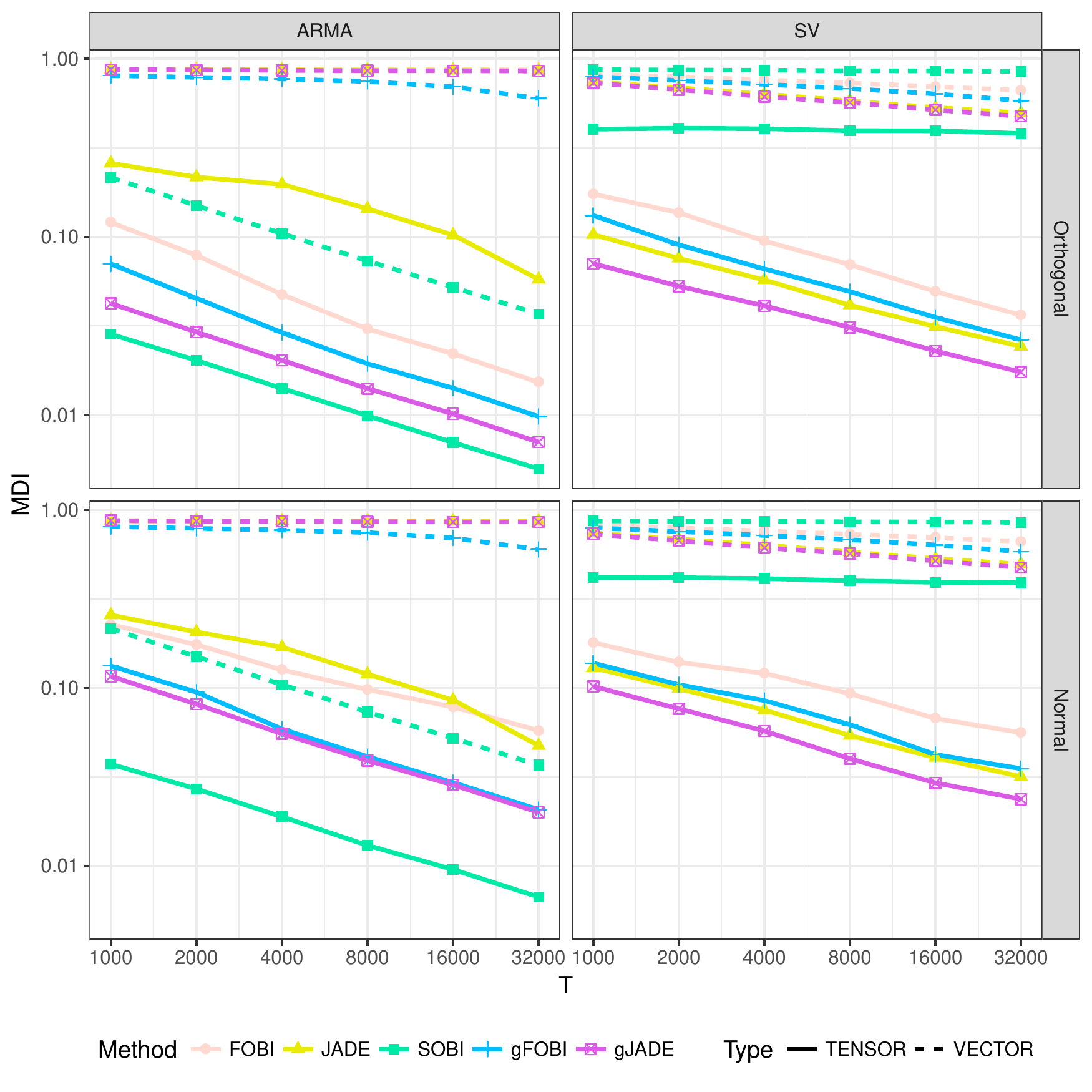}
\caption{The results of the efficiency comparison in the form of average MDI-values based on 1000 repetitions. Both axes have logarithmic scales.}
\label{fig:simul}
\end{figure}

The mean minimum distance indices (MDI) over different combinations of mixing, method, setting and series length are depicted in Figure \ref{fig:simul}. Based on the plots, the overall trend is that the vector methods offer no competition to the tensorial ones, the only exception being SOBI which in the ARMA setting outmatches TJADE (and also TFOBI in the case of normal mixing). Furthermore, in the ARMA setting the methods exploiting temporal dependence surpass the ones assuming the independence of observations and in the stochastic volatility setting the best performance is by a clear margin given by a method utilizing temporal dependence, TgJADE. Otherwise the plots offer no surprises: TSOBI, operating solely on autocovariances, outperforms the other methods in the ARMA setting and the FOBI and JADE-based tensor methods dominate in the stochastic volatility setting where no second-order information is present. In the latter case TJADE is almost as good as TgJADE indicating that most of the separation information is in the stochastic volatility setting contained in the higher marginal moments of the series and not in the temporal dependence. Of interest is also the better performance of gFOBI relative to the other non-SOBI vector methods in the ARMA setting, the same phenomenon that was observed also in \cite{MatilainenNordhausenOja2015}.
The figure also clearly shows that the tensor methods are not affine equivariant - but the dependence on the mixing matrix at least in this case seems to be negligible, especially when compared to the performance of the affine equivariant vectorized versions of the methods.

\subsection{Application: fMRI-data}

\begin{figure}[tp]
\centering
\includegraphics[width=1.0\textwidth]{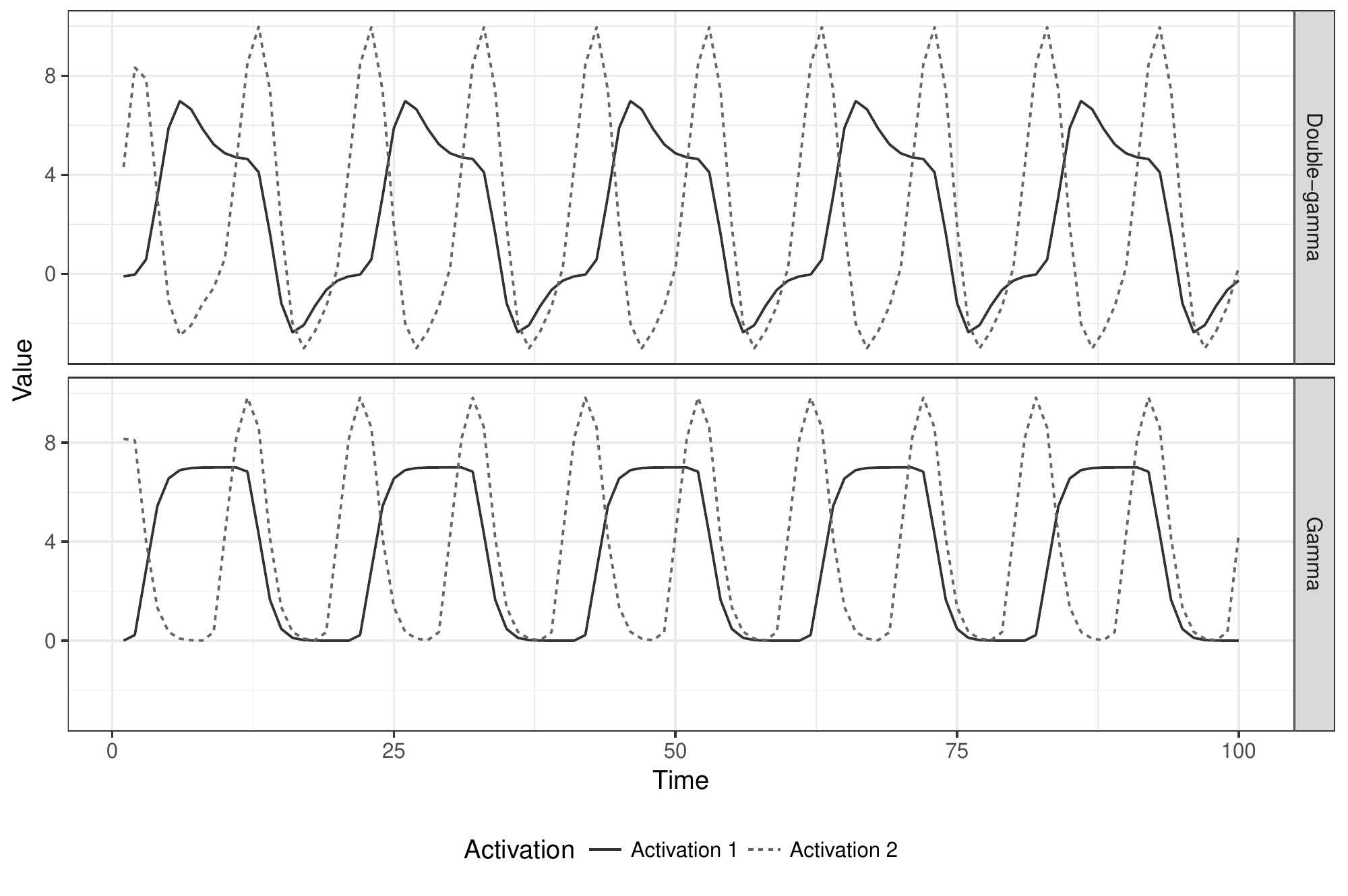}
\caption{The different activation signals we are trying to estimate from the fMRI-data. The upper plot corresponds to the double-gamma HRF and the lower to the gamma HRF.}
\label{fig:fmri_signals}
\end{figure}

We next consider an application of the methods to task-based fMRI-data simulated using the R-package \textit{neuRosim} \citep{RneuRosim}. The package offers a variety of options for simulating 3-dimensional, noisy, single-subject fMRI-measurements and as our setting we use a particular subset of the settings used in \cite{virta2016applying}. The simulated tensor-valued time series are of length $T = 100$ and the individual 3D-frames are of size $64 \times 64 \times 64$ containing two activation signals residing in disjoint regions, the estimation of which is our objective. The activation signals shown in Figure~\ref{fig:fmri_signals} are convolutions of a stimulus function and a \textit{haemodynamic response function} (HRF), with two different choices for the latter, gamma and double-gamma. See the section IV.A of \cite{virta2016applying} for more elaborate description of the simulation settings, the simulation code and images of the activation regions.

Our objective is then to estimate the activation signals in Figure \ref{fig:fmri_signals} and for this we compared a total of 8 methods: TFOBI, TJADE, TSOBI, TgFOBI, TgJADE and the vector-valued versions of all non-JADE-based methods, for which the data cubes were first vectorized. The reason for excluding JADE and gJADE is that their computation was not feasible on a computer with 24 cores and 64GB memory in reasonable time. As is standard in the processing of fMRI-data, prior to applying the methods we first reduced the dimensions of the observed series using principal component-based methods. For the tensor-valued methods \textit{tensorial PCA} \citep{virta2016applying} was used to reduce to the case of $6 \times 6 \times 2$ tensors and for the vector-valued methods SVD was used to compress the observations into vectors of length 72, which are then of comparable dimensionality with the reduced tensors.

\begin{figure}[tp]
\centering
\includegraphics[width=1.0\textwidth]{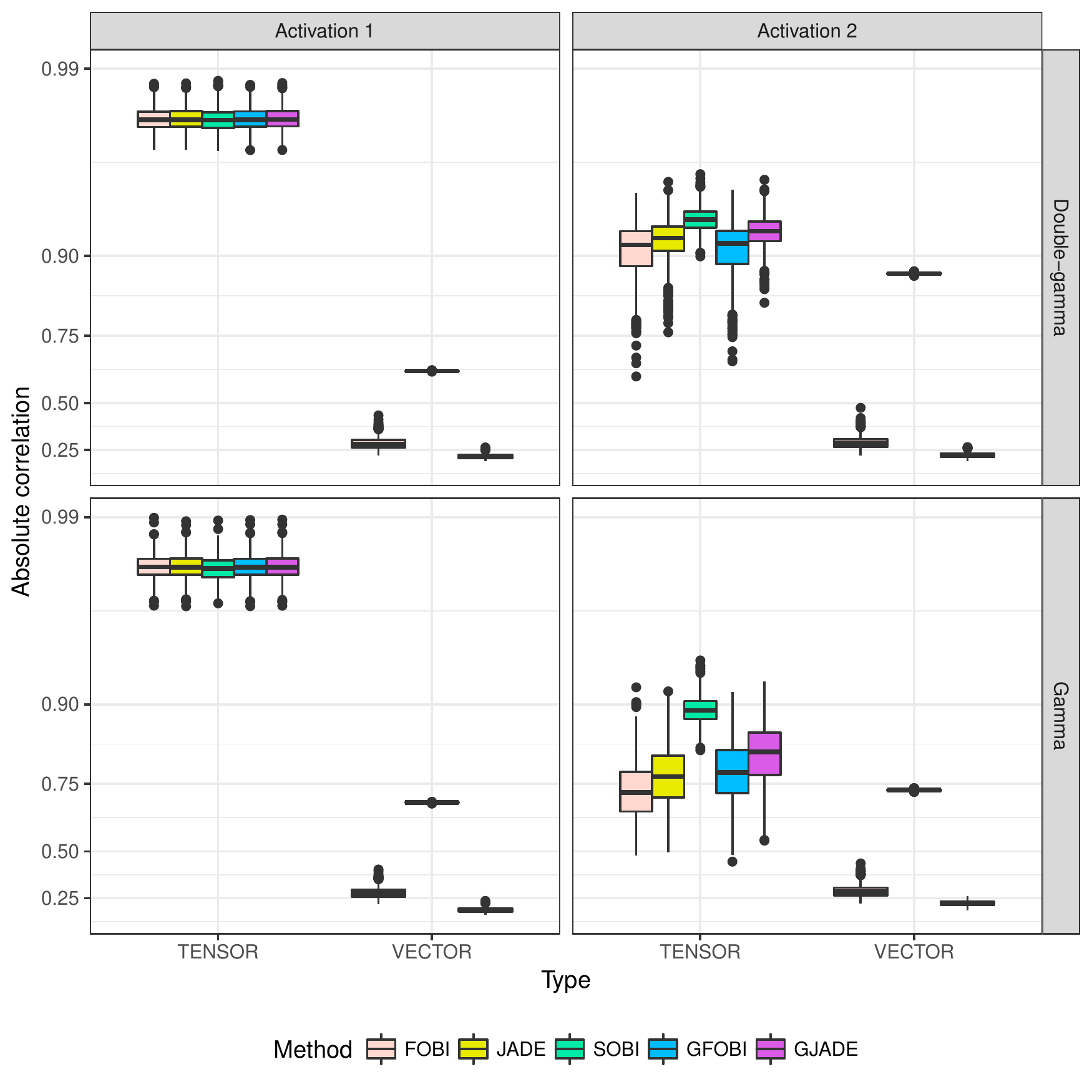}
\caption{The boxplots of resulting maximum absolute correlations between the activation signals and estimated components. In each of the four plots the boxplots correspond from left to right to TFOBI, TJADE, TSOBI, TgFOBI, TgJADE, FOBI, SOBI and gFOBI, respectively.}
\label{fig:fmri_results}
\end{figure}

We did 1000 replications under both described HRF-settings and for each replication recorded the highest absolute correlations with each of the true signals found amongst the obtained components. Thus the value of one indicates that the method succeeded exactly in extracting the corresponding signal. The resulting boxplots are shown in Figure \ref{fig:fmri_results}, with the y-axis obeying arctanh-scale (Fisher's z-transformation) to better emphasize the upper part of its range. The results clearly indicate that apart from SOBI the vectorial methods offer no competition to the tensorial methods. And while the first activation signal is for both HRFs estimated almost equally well by all tensor methods their performances differ in finding the second signal, TSOBI triumphing over the other methods, which is not too surprising as the signals of interest have rather regular patterns and do not really exhibit features of stochastic volatility. Note that the results for TFOBI, TJADE and FOBI in Figure~\ref{fig:fmri_results} match those obtained in the same settings in \cite{virta2016applying}.

\subsection{Application: video processing}
\begin{figure}[tp]
\centering
\includegraphics[width=1.0\textwidth]{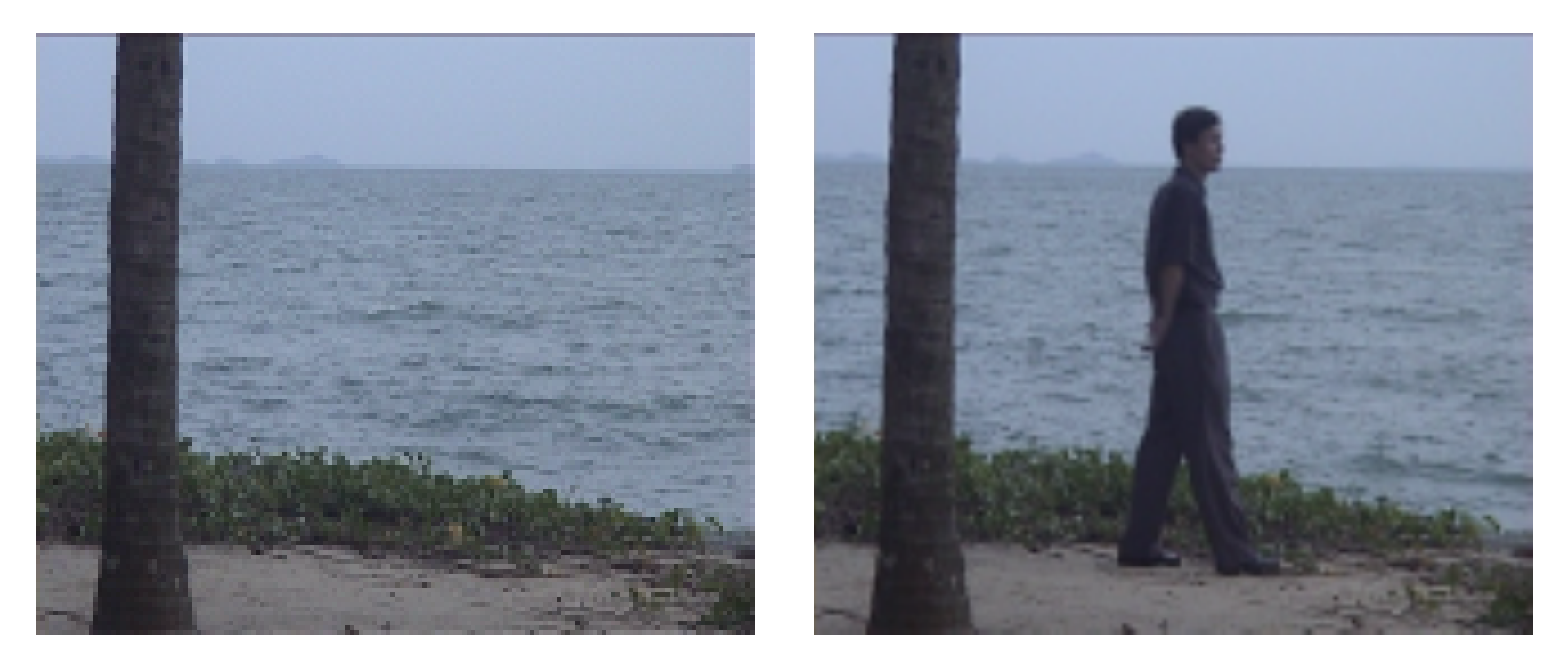}
\caption{The frames 1 and 520 of the original video on the left and right, respectively.}
\label{fig:video_frames}
\end{figure}

As described in the introduction, one of the more natural examples of tensor-valued time series is a video clip. A color video with resolution $h \times w$, 3 color channels and $T$ frames can be thought of as the sequence of realizations of a tensor-valued time series $\mathbb{X}_t \in \mathbb{R}^{h \times w \times 3}$ for $t = 1, \ldots , T$. As our third example we consider the WaterSurface video available online at \url{http://pages.cs.wisc.edu/~jiaxu/projects/gosus/supplement} that has been used for background subtraction in \cite{toyama1999wallflower, li2004statistical}. The same video was also used for outlier detection in \cite{rousseeuw2016measure} and they have the preprocessed video available at \url{https://wis.kuleuven.be/stat/robust/software}. The video data tensor consists of a total of $T = 633$ color frames of size $128 \times 160$ depicting a beach, sea and a tree, with a man entering the scene from left and passing the tree during frames 480-500 and staying in the picture for the rest of the clip. As an illustration Figure~\ref{fig:video_frames} shows frames 1 and 520 representing the typical view of the video without and with the man.

\begin{figure}[tp]
\centering
\includegraphics[width=0.95\textwidth]{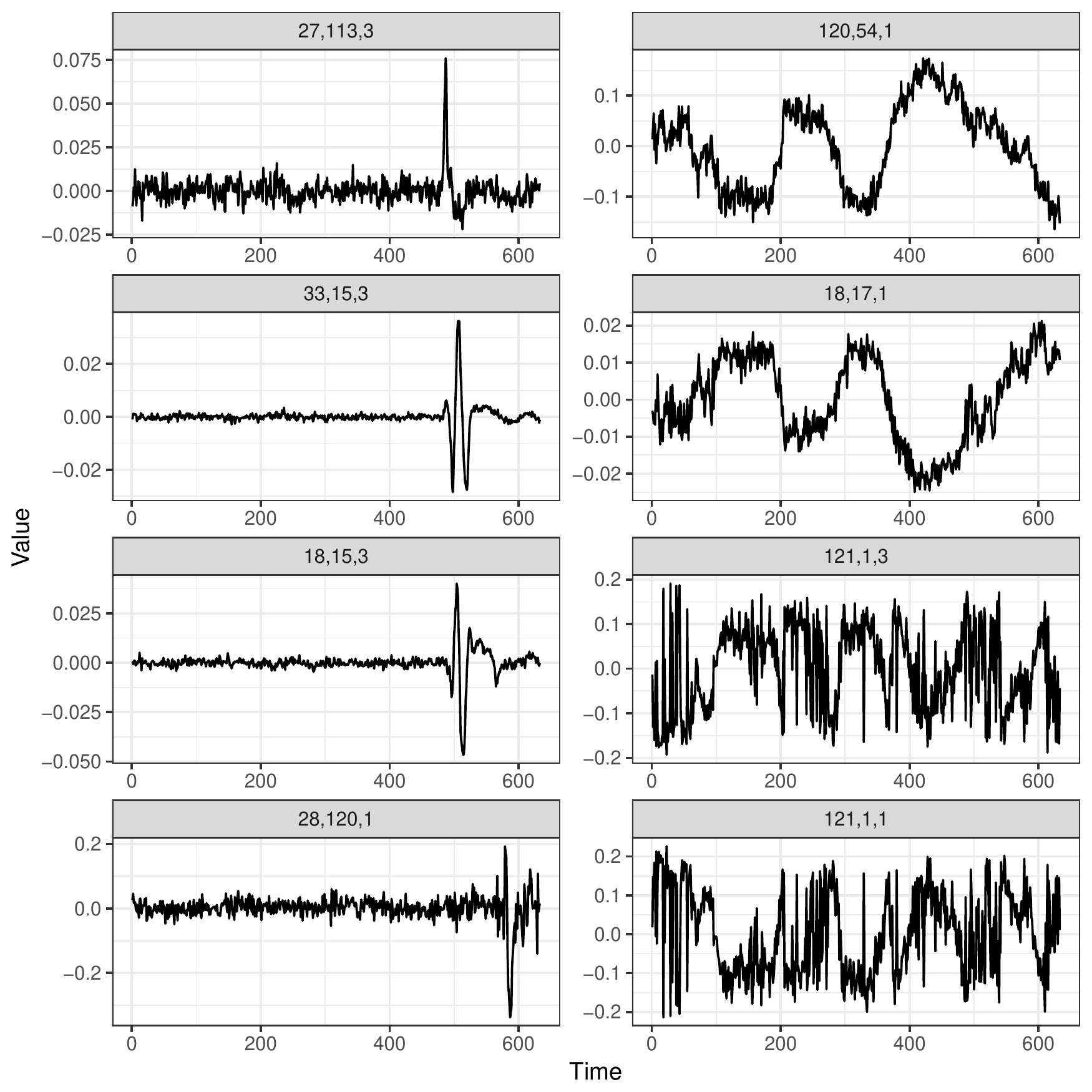}
\caption{On the left-hand side from top to bottom the four source series with the highest kurtosis values and on the right-hand side from bottom to top the four source series with the lowest kurtosis values. The titles of the source series refer to the indices of the series in the extracted source tensor $\mathbb{Z}_t$.}
\label{fig:video_sources}
\end{figure}

Our objective with the video is similar as in \cite{rousseeuw2016measure}, trying to detect the time window during which the man is in the scene. Assuming that the model~(\ref{eq:t_bssmodel}) holds in this case, the relative staticity of the video before the man's entrance makes it reasonable to assume that the change point is also visible in the extracted source signals in $\mathbb{Z}_t$. For demonstration purposes we considered for this example only TSOBI and extracted the four components with the highest kurtoses and the four components with the lowest kurtoses shown in Figure~\ref{fig:video_sources}. We chose kurtosis as a criterion for selecting the most interesting components as bimodal distributions are often characterized by extremal kurtosis values. The change point is now visible in the first three source series on the left-hand side of Figure \ref{fig:video_sources}, the latter two series having their peaks only after the man has already passed the tree. The fourth series on the left corresponds to a moment when the man has reached the edge of the screen and starts to slowly turn around. Note that contrary to the outlier detection method in \cite{rousseeuw2016measure} which flagged all frames after the change point as outlying, TSOBI produced signals that identify the frames during which big changes occur. Interestingly, the first two source series on the right also capture some systematic component in the video although exactly which is not evident from the video. Finally, the running time of the TSOBI-algorithm for the considered video tensor of size $128 \times 160 \times 3 \times 633$ is roughly ten minutes, making it a suitable method also for large-scale problems. Note that using vector-valued BSS methods without resorting to some form of \textit{a priori} dimension reduction is not possible in this case, the vectorized data set having far more variables ($p = 61440$) than observations ($T = 633)$.

\section{Discussion} \label{sec:disc}

Higher order structures are nowadays an increasingly popular form of data and manifest, for example, in the form of video data. One way to treat such observations is as a time series where at each time point we observe a tensor of the same size which, although seemingly natural, as a viewpoint is still largely unexplored in the literature. However, with the increasingly high-dimensional structures we also face the problem of separating the actual information content from the noise. In multivariate time series analysis a standard solution is blind source separation, a technique which in this paper was extended to tensor-valued time series.

Our approach is model-based and we first defined the tensor blind source separation model, a natural extension of the corresponding vector-valued model, where the objective is to recover a latent tensor-valued time series subjected to multimodal linear transformation. Next, we presented tensorial versions of three existing blind source separation methods, SOBI, gFOBI and gJADE, and proved that each of them is capable of estimating the latent source series under specific assumptions. The methods operate respectively on second-order, fourth-order and joint lagged fourth-order moments making the SOBI-based method specifically suited for the separation of linear processes and the last two methods suited to cases where no second-order information exists.

Simulation studies showed that under the model the proposed methods were highly superior not only to the often used combination of vectorizing the tensors and using standard multivariate methods but also to two tensorial dimension reduction methods designed for i.i.d data, TFOBI and TJADE. Hence, if the data exhibits a natural tensor structure and temporal correlation \textit{both} should also be taken into account in the analysis and processing of the data. In addition, two experiments conducted respectively on simulated fMRI-data and real video data revealed that the proposed methods are also adept in condensing information and extracting signals of importance from high-dimensional tensor-valued time series.

To also theoretically evaluate the gain in performance, our future plans include the investigation of whether asymptotical results like those obtained for SOBI in \citet{MiettinenIllnerNordhausenOjaTaskinenTheis2015} can be obtained also for the proposed tensor methods, especially TSOBI. Furthermore, we also plan to extend standard time series BSS methods that allow relaxing the stationarity assumption, such as \citet{ChoiCichocki2000a,ChoiCichocki2000b,Nordhausen2014}, to the tensorial case. Another interesting problem is the choice of the most important source components in applications. In Section \ref{sec:simu} we used absolute correlation and extremal kurtosis as criteria but several others, such as the number of change points, could be used as well. Also this we plan to explore further in our future work.

\section{Acknowledgements}

The authors would like to thank the two anonymous referees for their helpful comments and insightful suggestions.

\appendix

\section{The proofs and simulation distributions} \label{sec:appe}

\begin{proof}[The proof of Theorem \ref{theo:t_diag}]
Note first that equation (17) in \citet{VirtaLiNordhausenOja2016b} says that the $m$-flattening $(\bo{X}^{st}_t)^{(m)}$ is equal to $d \cdot \bo{U}_m \bo{Z}^{(m)}_t (\bo{U}_{m+1} \otimes \cdots \otimes \bo{U}_r \otimes \bo{U}_1 \otimes \cdots \otimes \bo{U}_{m-1})$, where $d$ is the constant of proportionality and the Kronecker product as a product of orthogonal matrices is itself an orthogonal matrix.

Theorem \ref{theo:t_rotation} then implies that $\mathbb{X}^{st}_t \propto \mathbb{Z}_t \odot_1 \bo{U}_1 \cdots \odot_r \bo{U}_r$ which, when plugged in into $\bo{\Sigma}_\tau^m (\mathbb{X}^{st}_t)$ and $\bo{B}_\tau^m (\mathbb{X}^{st}_t)$ in combination with the previous result, \ref{assu:t_tsobiworks} and \ref{assu:t_tgfobiworks} directly yields the result for the parts \textit{i)} and \textit{ii)}.

To prove part \textit{iii)}, it is straightforward to show that for any lags $\tau_1, \tau_2, \tau_3, \tau_4$ and indices $i, j = 1, \ldots ,p_m$ we have
\[
\bo{B}_{\tau_1 \tau_2 \tau_3 \tau_4 i j}^m (\mathbb{X}^{st}_t) = d^4 \sum_{k,l} u^{(m)}_{ik} u^{(m)}_{jl} \cdot \textbf{U}_m \bo{B}_{\tau_1 \tau_2 \tau_3 \tau_4 k l}^m (\mathbb{Z}_t) \bo{U}_m^T,
\]
where $(\bo{U}_m)_{kl} = u^{(m)}_{kl}$. By the assumption we have for any fixed $\tau, i, j$ that $\bo{C}^m_{\tau i j}(\mathbb{Z}_t) = \bo{D}^m_{\tau i j}$ is diagonal and the matrix $\bo{C}_{\tau i j}^m (\mathbb{X}^{st}_t)$ can thus be written as
\begin{align*}
\bo{C}_{\tau i j}^m (\mathbb{X}^{st}_t) &= d^4 \cdot \bo{U}_m \left( \sum_{k,l} u^{(m)}_{ik} u^{(m)}_{jl} \cdot \bo{D}^m_{\tau i j} \right. \\
&+ \bo{\Sigma}^m_0(\mathbb{Z}_t)
\left[ \sum_{k,l} u^{(m)}_{ik} u^{(m)}_{jl} \cdot (\bo{E}^{kl} + \bo{E}^{lk} + \bo{I} ) \right. \\
&- \left. \left. \bo{U}_m^T (\bo{E}^{ij} + \bo{E}^{ji} + \bo{I} ) \bo{U}_m \right] \bo{\Sigma}^m_0(\mathbb{Z}_t) \right) \bo{U}^T_m,
\end{align*}
where $\bo{\Sigma}^m_0(\mathbb{Z}_t)$ is diagonal and by element-wise inspection the quantity inside the square brackets can be showed to be diagonal as well. Thus the right-hand side expression is the eigendecomposition of $\bo{C}_{\tau i j}^m (\mathbb{X}^{st}_t)$ proving the claim.
\end{proof}

\begin{proof}[The models used in the simulation study]\renewcommand{\qedsymbol}{}
In the ARMA setting the 12 used time series models were: AR(1) with $\phi = 0.9$; AR(1) with $\phi = -0.9$; MA(2) with $\bo{\theta} = (0.5, -0.5)$; AR(2) with $\bo{\phi} = (-0.5, -0.3)$; ARMA(4,2) with $\bo{\phi} = (0.5, -0.3, 0.1, -0.1)$ and $\bo{\theta} = (0.7, -0.3)$; ARMA(2,4) with $\bo{\phi} = (-0.7, 0.1)$ and $\bo{\theta} = (0.9, 0.3, 0.1, -0.1)$ and MA(5), MA(10), MA(20), MA(30), MA(40) and MA(50) with parameters sampled independently from the uniform distribution on $(-1, 1)$.

In the SV setting the used models consisted of both GARCH and stochastic volatility models, SV($\mu, \phi, \sigma, \nu$), where $\mu$, $\phi$ and $\sigma$ are respectively the mean, AR-parameter and the volatility of the latent AR(1)-process and $\nu$ is the number of degrees of freedom of the innovation $t$-distribution, see \cite{kastner2014ancillarity}. The six used stochastic volatility models were then SV($-10, 0.98, 0.2, \infty$), SV($-5, -0.98, 0.2, 10$), SV($-10, 0.7, 0.7, \infty$), SV($-5, -0.70, 0.7, 10$), SV($-9, 0.20, 0.01, \infty$) and SV($-9, -0.20, 0.01, 10$). In addition, the six used GARCH models were: ARCH(1) with $\alpha = 0.7$; GARCH(1,1) with $\alpha = 0.2$ and $\beta = 0.2$; GARCH(1,1) with $\alpha = 0.1$ and $\beta = 0.8$; ARCH(4) with $\bo{\alpha} = (0.20, 0.10, 0.05, 0.01)$; GARCH(3,1) with $\bo{\alpha}=(0.05, 0.03, 0.01)$ and $\beta=0.5$ and ARCH(10) with $\bo{\alpha} = (0.20, 0.14, 0.12, 0.10, 0.05, 0.05, 0.04, 0.03, 0.02, 0.01)$, see \cite{bollerslev1986generalized} for more details on ARCH and GARCH models.
\end{proof}

\section*{References}

\bibliography{new_references}

\end{document}